\documentclass[12pt]{article}
\usepackage{amssymb}

\usepackage{amsfonts,amssymb,eucal,amsmath}
\title{ Notes on a Theorem of Benci-Gluck-Ziller-Hayashi}

\author{\small\sc Fengying Li\footnote{Email:lify0308@163.com} and Shiqing Zhang\footnote{Email:zhangshiqing@msn.com}\\
 \small \it  The School of Economic and Mathematics,
Southwestern
University of Finance and Economics, \\
\small \it Chengdu 611130, China\\
\small \it Mathematical School, Sichuan University, Chengdu610064
,China}

\date{}

\begin{document}
\maketitle

\begin{abstract}

We use constrained variational minimizing methods to study the existence of periodic solutions with a prescribed energy for
a class of second order Hamiltonian systems with a $C^2$ potential function which may have an unbounded potential well. Our
result can be regarded as complementary to the well-known
theorem of Benci-Gluck-Ziller and Hayashi.\\
\end{abstract}



\section*{1. Introduction}
\setcounter{section}{1} \setcounter{equation}{0}

Based on the earlier works of Seifert([20]) in 1948 and
Rabinowitz([18,19]) in 1978 and 1979, Benci ([4]), and
Gluck-Ziller ([11]), and Hayashi([13]) published work examining the
periodic solutions for second order Hamiltonian systems
\begin{equation}
\ddot{q}+V^{\prime}(q)=0\label{1.1}
\end{equation}
\begin{equation}
\frac{1}{2}|\dot{q}|^2+V(q)=h\label{1.2}
\end{equation}
with a fixed energy. Utilizing the Jacobi metric and very
complicated geodesic methods with algebraic topology, they proved the following general theorem:
\vspace{0.3cm}

{\bf Theorem 1.1}\ \ Suppose $V\in C^2(R^n,R)$. If the potential
well
$$\{x\in R^n|V(x)\leq h\}$$
 is bounded and non-empty,
then the system (1.1)-(1.2) has a periodic solution with energy h.

  Furthermore, if
$$V^{\prime}(x)\not=0, \hspace{0.4cm} \forall x\in\{x\in R^n|V(x)=h\},$$
then the system (1.1)-(1.2) has a nonconstant periodic solution with energy h.
\vspace{0.4cm}

 For the existence of multiple periodic solutions for (1.1)-(1.2)
 with compact energy surfaces, we can refer to Groessen([12]) and Long[14]
 and the references therein.

 In 1987, Ambrosetti-Coti Zelati[2] successfully used Clark-Ekeland's dual action principle and Ambrosetti-Rabinowitz's Mountain Pass theorem to study the existence of $T$-periodic solutions of
 the second-order equation
 $$-\ddot x=\nabla U(x),$$
where $$U=V\in C^2(\Omega;\bold R)$$
such that
 $$U(x)\to\infty,x\to \Gamma=\partial\Omega;$$
with  $\Omega\subset\bold R^n$ a bounded convex domain. Their principle result is the following:
\vspace{0.3cm}

{\bf Theorem 1.2}\ \ Suppose

\begin{enumerate}

    \item  $U(0)=0=\min U$

    \item $U(x)\le \theta(x,\nabla U(x))$ for some $\theta\in (0,\tfrac 12)$ and for all $x$ near $\Gamma$
   (superquadraticity near $\Gamma$)

   \item $(U''(x)y,y) \ge k|y|^2$ for some $k>0$ and for
   all $(x,y)\in\Omega\times \bold R^N$.
\end{enumerate}

    Let $\omega_N$ be the greatest eigenvalue  of $U'' (0)$ and $T_0=(2/\omega_N)^{1/2}$.

    Then  $-\ddot x=\nabla U(x)$ has for each $T\in (0,T_0)$ a periodic solution with minimal period $T$.
\vspace{0.4cm}

The dual variational principle and Mountain Pass Lemma again proved the essential ingredients for the following theorem of Coti
Zelati-Ekeland-Lions [8] concerning Hamiltonian systems in convex
potential wells.
\vspace{0.3cm}

{\bf Theorem 1.3}\ \ Let $\Omega$ be a convex open subset of $R^n$
containing the origin $O$. Let $V\in C^2(\Omega,R)$ be such that

$(V1).\ V(q)\geq V(O)=0,\forall q\in\Omega$

$(V2).\ \forall q\neq O,V''(q)>0$

$(V3).\ \exists \omega>0,$ such that

$$V(q)\leq \frac{\omega}{2}\|q\|^2,\forall \|q\|<\epsilon$$

and
\vspace{0.1cm}

$(V4). V''(q)^{-1}\rightarrow 0,\|q\|\rightarrow 0 $
 or,

$(V4)' .V''(q)^{-1}\rightarrow 0, q\rightarrow\partial\Omega.$

Then, for every $T<\frac{2\pi}{\sqrt{\omega}}$, the system (1.1) has a
solution with minimal period $T.$
\vspace{0.4cm}

 In Theorems 1.2 and 1.3, the authors assumed the convex conditions
 for potentials and potential wells in order to apply Clark-Ekeland's
 dual variational principle. We observe that
  Theorems 1.1-1.3 essentially make the
 assumption
 $$V(x)\to\infty,x\to \Gamma=\partial\Omega$$
 so that all potential wells are bounded.
 We wish to generalize Theorems 1.1-1.3 from two directions:
 (1) We dispense with the convex assumption on potential functions,
 (2) $V(x)$ can be uniformly bounded, and the potential well can be
 unbounded.

 In 1987, D.Offin ([16]) generalized Theorem 1.1 to some
    non-compact cases for $V\in C^3(R^n,R)$ under complicated geometric
    assumptions on the potential wells; however, these geometric conditions
    appear difficult to verify for concrete potentials.
In 2009, Berg-Pasquotto-Vandervorst ([5]) studied the closed
orbits on non-compact manifolds with some complex topological
assumptions.

Using simpler constrained variational minimizing method, we obtain the following result:
\vspace{0.2cm}

{\bf Theorem 1.4}\ \ Suppose $V\in C^2(R^n,R),h\in R$ satisfies

$(V_1).\  V(-q)=V(q)$

$(V_2).\ V^{\prime}(q)q>0,\forall q\neq 0$

$(V_3).\  3V^{\prime}(q)q+(V''(q)q,q)\neq 0,\forall q\neq 0$

 $(V_4).\ \exists
\mu_1>0,\mu_2\geq 0,$ such that

$$V^{\prime}(q)\cdot q\geq\mu_1 V(q)-\mu_2$$

$(V_5).\
\lim_{|q|\rightarrow\infty}Sup[V(q)+\frac{1}{2}V^{\prime}(q)q]\leq
A$

$(V_6).\ \frac{\mu_2}{\mu_1}<h<A.$

 Then the system $(1.1)-(1.2)$ has at least
one non-constant periodic solution with the given energy h.\\

{\bf Corollary 1.5}\ \ Suppose $V(q)=a|q|^{2n},a>0$, then the system $\forall
h>0$,$(1.1)-(1.2)$ has at least
one non-constant periodic solution with the given energy h.\\

{\bf Remark 1}\ \ Suppose $V(x)$ is the following well-known
$C^{\infty}$ function:

$$V(x)=e^{\frac{-1}{|x|}},\forall x\neq 0;$$
$$V(0)=0.$$

Then $V(x)$ satisfies $(V_1)-(V_5)$ if we take $\mu_1=\mu_2>0$ and
$A=1$ in Theorem 1.4, but $(V_6)$ does not hold .\\

{\bf Proof}  In fact,it's easy to check $(V_1)-(V_5)$:

(1). It's obvious for $(V_1)$.

(2). For $(V_2)$ and $(V_3)$, we notice that

$$V^{\prime}(x)x=\frac{1}{|x|}e^{\frac{-1}{|x|}}>0,\forall x\not=0,$$

$$(V''(x)x,x)=e^{\frac{-1}{|x|}}(\frac{-2}{|x|}+\frac{1}{|x|^2})$$

$$3V^{\prime}(x)x+(V''(x)x,x)=e^{\frac{-1}{|x|}}(\frac{1}{|x|}+\frac{1}{|x|^2})
>0,\forall x\neq 0.$$

(3). For $(V_4)$, we set

$$w(x)=(\frac{1}{|x|}-\mu_1)e^{\frac{-1}{|x|}}; \hspace{0.2cm} x\not=0, w(0)=0.$$

We will prove $w(x)>-\mu_1$; in fact,

$$w^{\prime}(x)=[\frac{1}{|x|}-(1+\mu_1)]\frac{x}{|x|^3}e^{\frac{-1}{|x|}}; x\not=0, w^{\prime}(0)=0.$$
From $w^{\prime}(x)=0$ ,we have $x=-\frac{1}{1+\mu_1}$ or $0$ or
$\frac{1}{1+\mu_1}$.

 It's easy to see that $w(x)$ is strictly increasing on $(-\infty,-\frac{1}{1+\mu_1}]$ and
 $[0,\frac{1}{1+\mu_1}]$ but  strictly decreasing on $[\frac{-1}{1+\mu_1},0]$ and $[\frac{1}{1+\mu_1},+\infty)$.
 We notice that
 $$\lim_{|x|\rightarrow +\infty}w(x)=-\mu_1,$$
 and
 $$w(0)=0.$$
 So
$$w(x)>-\mu_1.$$
 When we take $\mu_2=\mu_1>0$,$(V_4)$ holds.

(4). For $(V_5)$, we have

$$V(x)+\frac{1}{2}V^{\prime}(x)x=e^{\frac{-1}{|x|}}(1+\frac{1}{2}\frac{1}{|x|})
< 1,\forall x\neq 0;$$

$$ V(0)+\frac{1}{2}V^{\prime}(0)0=0.$$

{\bf Corollary 1.6}\ \ Given any $a>0,n\in N$, suppose
$V(x)=a|x|^{2n}+e^{\frac{-1}{|x|}}, x\not=0, V(0)=0$. Then $\forall
h>1$, the system $(1.1)-(1.2)$ has at least
one non-constant periodic solution with the given energy $h$.\\

 {\bf Remark 2}\ \ The potential $V(x)$ in Remark 1 is noteworthy since the
 potential function is non-convex and bounded which satisfies neither of
 the conditions of Theorems 1.1-1.3, Offin's geometrical
 conditions, nor Berg-Pasquotto-Vandervorst's complex topological
assumptions.
 Notice the special properties for our potential well. It is
 a bounded set if $h<1$, but for $h\geq 1$ it is $R^n$ - an unbounded set.
 We also notice that the symmetrical
 condition on the potential simplified our Theorem 1.4 and it's
 proof; it seems interesting to observe to obtain
 non-constant periodic solutions if the symmetrical condition is
 deleted.

\section{ A Few  Lemmas}
\setcounter{section}{2} \setcounter{equation}{0}

 Let
 $$H^1=W^{1,2}(R/Z,R^n)=\{u:R\rightarrow R^n,u\in L^2,\dot{u}\in L^2,u(t+1)=u(t)\}$$
 Then the standard $H^1$ norm is equivalent to
 $$\|u\|=\|u\|_{H^1}=\left(\int^1_0|\dot{u}|^2dt\right)^{1/2}+|\int_0^1u(t)dt|.$$

{\bf Lemma 2.1}([1])\ \ Let
$$M=\{u\in H^1|\int_0^1(V(u)+\frac{1}{2}V^{\prime}(u)u)dt=h\}.$$
If $(V_3)$ holds, then $M$ is a $C^1$ manifold with codimension 1
in $H^1.$

Let

$$f(u)=\frac{1}{4}\int^1_0|\dot{u}|^2dt\int^1_0V^{\prime}(u)udt$$

and $\widetilde{u}\in M$ be such that
$f^{\prime}(\widetilde{u})=0$
 and $f(\widetilde{u})>0$. Set

$$\frac{1}{T^2}=\frac{\int^1_0V^{\prime}(\widetilde{u})\widetilde{u}dt}{\int^1_0|\dot{\widetilde{u}}|^2dt}$$
If $(V_2)$ holds, then $\widetilde{q}(t)=\widetilde{u}(t/T)$ is a
non-constant $T$-periodic solution for (1.1)-(1.2).
\vspace{0.4cm}

When the potential is even, then by Palais's symmetrical principle
([17]) and Lemma 2.1, we have
\vspace{0.2cm}

{\bf Lemma 2.2}([1])\ \ Let
$$F=\{u\in M|u(t+1/2)=-u(t)\}$$
and suppose $(V_1)-(V_3)$ holds. If $\widetilde{u}\in F$ be such that
$f^{\prime}(\widetilde{u})=0$ and $f(\widetilde{u})>0$,then
$\widetilde{q}(t)=\widetilde{u}(t/T)$ is a non-constant
$T$-periodic solution for (1.1)-(1.2). In addition, we have
$$\forall u\in F,\int_0^1u(t)dt=O.$$
\vspace{0.1cm}

Recall the following two classic results.
\vspace{0.3cm}

 {\bf Lemma
2.3}(Sobolev-Rellich-Kondrachov[15],[22])
$$W^{1,2}(R/Z,R^n)\subset C(R/Z,R^n)$$
and the imbedding is compact.\\
\vspace{0.2cm}

{\bf Lemma 2.4}(Eberlein-Smulian [21])\ \ A Banach space $X$ is
reflexive if and only if any bounded sequence in $X$ has a weakly
convergent subsequence.\\

{\bf Definition 2.1}(Tonelli ,[15])\ \ Let $X$ is a Banach
space and $M\subset X$. If it the case that for any sequence $\{x_n\}\subset M$ strongly
convergent to $x_0$ ($x_n \rightarrow x_0$), we have $x_0\in M$, then
we call $M$ a strongly closed (closed) subset of $X$; if for
any $\{x_n\}\subset M$ weakly convergent to
$x_0$ ($x_n\rightharpoonup x_0$), we have $x_0\in M$, then we call
$M$ a weakly closed subset of $X$.

Let $f:M\rightarrow R$.

(i). If for any $\{x_n\}\subset M$ strongly convergent to
$x_0$,we have

 $$liminf f(x_n)\geq f(x_0),$$
 then we say $f(x)$ is lower  semi-continuous at $x_0$.

(ii). If for any $\{x_n\}\subset M$  weakly convergent to
$x_0$, we have

 $$liminf f(x_n)\geq f(x_0),$$

 then we say $f(x)$ is weakly lower semi-continuous at $x_0$.
\vspace{0.3cm}

 Using his variational principle, Ekeland proved
\vspace{0.3cm}

{\bf Lemma 2.5}(Ekeland[9])\ \ Let $X$ be a Banach
space and $F\subset X$ a closed (weakly closed) subset. Suppose that $\Phi $ defined on $X$ is Gateaux-differentiable and  lower
 semi-continuous (or weakly lower semi-continuous)
  and that $\Phi|_F$ restricted on $F$ is bounded from
 below. Then there is a sequence $x_n\subset F$ such that
 $$\Phi(x_n)\rightarrow\inf_{F}\Phi \hspace{0.4cm} \hbox{ and } \hspace{0.4cm} \|\Phi|_F^{'}(x_n)\|\rightarrow 0.$$

  {\bf Definition 2.2}([9,10])\ \ Let $X$ be a  Banach
space and $F\subset X$ a closed (weakly closed) subset.
 Suppose that $\Phi $ defined on $X$ is Gateaux-differentiable.
 If it is true that whenever $\{x_n\}\subset F$ such that

 $$\Phi(x_n)\rightarrow c \hspace{0.4cm} \hbox{ and } \hspace{0.4cm} \|\Phi|_F^{'}(x_n)\|\rightarrow 0,$$
then $\{x_n\}$ has a strongly convergent (weakly convergent) subsequence, we say $\Phi$ satisfies the $(PS)_{c,F}$ ($(WPS)_{c,F}$) condition at the level $c$
for the closed subset $F\subset X$.\\
\vspace{0.1cm}

Using $\bf Lemma2.5$, it is easy to prove the following lemma.
\vspace{0.3cm}

{\bf Lemma 2.6}\ \ Let $X$ be a  Banach space,

(i). Let $F\subset X$ be a closed  subset.
 Suppose that $\Phi $ defined on $X$ is Gateaux-differentiable and  lower
 semi-continuous and bounded from below on $F$.
 If $\Phi$ satisfies $(PS)_{\inf\Phi,F}$ condition, then $\Phi$ attains its infimum on $F$.

 (ii).Let $F\subset X$ be a weakly closed subset. Suppose that $\Phi $ defined on $F$ is
 Gateaux-differentiable and weakly lower semi-continuous
 and bounded from below on $F$.
 If $\Phi$ satisfies $(WPS)_{\inf\Phi,F}$ condition, then $\Phi$ attains its infimum on $F$.

\section{The Proof of Theorem 1.4}
\setcounter{section}{3} \setcounter{equation}{0}

We prove the Theorem as a sequence of claims.
\vspace{0.2cm}

 {\bf Claim 3.1} If $(V_1)-(V_6)$ hold, then for any given $c>0$,
 $f(u)$ satisfies the $(PS)_{c,F}$
 condition; that is, if $\{u_n\}\subset F$ satisfies
\begin{eqnarray}
 f(u_n)\rightarrow c>0 \ \ \hbox{ and }\ \
f|_F^{\prime}(u_n)\rightarrow 0,\label{3.1}
\end{eqnarray}
then $\{u_n\}$ has a strongly convergent subsequence.\

{\bf Proof}\ \  First, we prove the constrained set
$F\not=\emptyset$ under our assumptions. Using the notation
of [1], for $a>0$ let

\begin{eqnarray}
g_u(a)=g(au)=\int^1_0[V(au)+\frac{1}{2}V^{\prime}(au)au] dt.
\label{3.2}
\end{eqnarray}
By the assumption $(V_3)$, we have

\begin{eqnarray}
\frac{d}{da}g_u(a)\not=0
\label{3.2}
\end{eqnarray}
and so $g_u$ is strictly monotone. By $(V_5)$, we have

\begin{eqnarray}
\lim_{a\rightarrow +\infty}g_u(a)\leq A
 \label{3.2}
\end{eqnarray}
By $(V_4)$, we notice that
\begin{eqnarray}
g_u(0)=V(O)\leq\frac{\mu_2}{\mu_1}.
 \label{3.2}
\end{eqnarray}
So for $V(O)<h<A$, the equation $g_u(a)=h$ has a unique solution
$a(u)$ with $a(u)u\in M.$

By $f(u_n)\rightarrow c$, we have
\begin{eqnarray}
\frac{1}{4}\int^1_0|\dot{u_n}(t)|^2dt\cdot\int^1_0V^{\prime}(u_n)u_ndt\rightarrow
c,\label{3.2}
\end{eqnarray}

and by $(V_4)$ we have

\begin{eqnarray}
h=\int^1_0(V(u_n)+\frac{1}{2}<V^{\prime}(u_n),u_n>)dt\leq
(\frac{1}{\mu_1}+\frac{1}{2})\int_0^1V^{\prime}(u_n)u_ndt+\frac{\mu_2}{\mu_1}\label{3.3}.
\end{eqnarray}

By (3.6) and (3.7) we have

\begin{eqnarray}
\int_0^1V^{\prime}(u_n)u_ndt\geq
\frac{h-\frac{\mu_2}{\mu_1}}{\frac{1}{2}+\frac{1}{\mu_1}}.\label{3.4}
\end{eqnarray}

Condition $(V_6)$ provides $h>\frac{\mu_2}{\mu_1}$. Then (3.6) and (3.8) imply
$\int^1_0|\dot{u_n}(t)|^2dt$ is bounded and
$\|u_n\|=\|\dot{u}_n\|_{L^2}$ is bounded.

We know that $H^1$ is a
reflexive Banach space, so by the embedding theorem, $\{u_n\}$ has a weakly convergent
subsequence which uniformly strongly
converges to $u\in H^{1}$. The argument to show $\{u_n\}$ has a
strongly convergent subsequence is standard, and we can refer to
Lemma 3.5 of Ambrosetti-Coti Zelati [1].
\vspace{0.2cm}

 {\bf Claim 3.2} $f(u)$ is weakly lower semi-continuous on
   $F$.

{\bf Proof} For any $u_n\subset F$ with $u_n\rightharpoonup u$, by
Sobolev's embedding Theorem we have the uniform convergence:

$$|u_n(t)-u(t)|_{\infty}\rightarrow 0.$$
Since $V\in C^1(R^n,R)$, we have
$$|V(u_n(t))-V(u(t))|_{\infty}\rightarrow 0.$$
By the weakly lower semi-continuity of norm, we have

$$\liminf (\int^1_0|\dot{u}_n|^2dt)^{\frac{1}{2}}\geq (\int^1_0|\dot{u}|^2dt)^{\frac{1}{2}} .$$
Calculating we see
 $$\liminf (\int^1_0|\dot{u}_n|^2dt)\geq\int^1_0|\dot{u}|^2dt,$$
and
$$\liminf f(u_n)=\liminf\frac{1}{4}\int^1_0|\dot{u_n}|^2dt\int^1_0V^{\prime}(u_n)u_ndt$$
$$\geq\frac{1}{4}\int^1_0|\dot{u}|^2dt\int^1_0V^{\prime}(u)udt=f(u).$$
\vspace{0.1cm}

 {\bf Claim 3.3} $F$ is a  weakly closed subset in $H^1$.

 {\bf Proof} This follows easily from Sobolev's embedding Theorem and $V\in C^1(R^n,R)$.
\vspace{0.2cm}

 {\bf Claim 3.4} The functional $f(u)$ has
positive lower bound on $F$

{\bf Proof} By the definitions of $f(u)$ and $F$ and the
assumption $(V_2)$, we have
$$f(u)=\frac{1}{4}\int^1_0|\dot{u}|^2dt\int^1_0(V^{\prime}(u)u)dt\geq 0,\forall u\in F.$$
Furthermore, we claim that
$$\inf f(u)>0 ;$$
otherwise,  $u(t)=const$, and by the symmetrical property
$u(t+1/2)=-u(t)$ we have $u(t)=0,\forall t\in R$. But by
assumptions $(V_4)$ and $(V_6)$ we have

$$ V(0)\leq\frac{\mu_2}{\mu_1}<h,$$
which contradicts the definition of $F$ since $V(0)=h $ if we have $0\in F$.
 Now by Lemmas 3.1-3.4 and Lemma 2.6, we see that $f(u)$
attains the infimum on $F$, and we know that the minimizer is
nonconstant.

\section*{Acknowledgements}

 The authors sincerely thank Professor P.Rabinowitz who brought the paper of D. Offin ([16]) to our attention.

\end{document}